\documentclass[12pt]{article}
\usepackage{amsxtra}
\usepackage{amssymb, amsmath}
\usepackage{hyperref}
\usepackage {graphicx}

\begin{document}
\begin{center}
\textbf{\Large{DIGITAL ERA: Universal Bimagic Squares}}
\end{center}

\bigskip
\begin{center}
\textbf{\large{Inder Jeet Taneja}}\\
Departamento de Matem\'{a}tica\\
Universidade Federal de Santa Catarina\\
88.040-900 Florian\'{o}polis, SC, Brazil.\\
\textit{e-mail: ijtaneja@gmail.com\\
http://www.mtm.ufsc.br/$\sim$taneja}
\end{center}

\bigskip
\begin{abstract}
\textit{In this short note we have produced different kinds of }\textbf{\textit{bimagic squares}}\textit{ using only the digits 0, 1, 2, 5 and 8. The universal bimagic squares presented are of order $8\times 8$, $9\times 9$, $16\times 16$ and $25\times 25$. In order to bring universal bimagic square of order $8\times 8$, we used only the digits 2 and 5. For the order $9\times 9$, we used only the digits 2, 5 and 8. For the universal bimagic square of order $16\times 16$ we used the digits 1, 2, 5 and 8 and finally for the order $25\times 25$, we used five digits 0, 1, 2, 5 and 8. In order to produce these universal bimagic squares we have used the digits in the digital form. }
\end{abstract}

\section{Introduction}

It is well known that there are digits specially used in watches, elevators,
etc. These digits are of type:
\begin{figure}[htbp]
\centerline{\includegraphics[width=1.63in,height=0.17in]{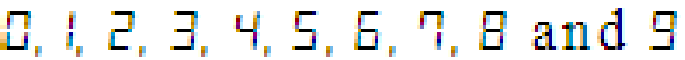}.}
\label{fig1}
\end{figure}

Some of these algarisms have special property such as
when we rotate them to $180^{o}$ degrees they become known algarisms again. These are
\begin{figure}[htbp]
\centerline{\includegraphics[width=1.16in,height=0.17in]{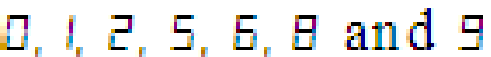}.}
\label{fig2}
\end{figure}

This means that the digits
\begin{figure}[htbp]
\centerline{\includegraphics[width=0.87in,height=0.17in]{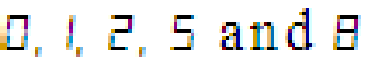}.}
\label{fig3}
\end{figure}

\noindent
remains the same while
\begin{figure}[htbp]
\centerline{\includegraphics[width=1.92in,height=0.17in]{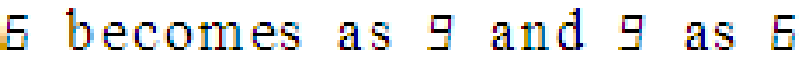}.}
\label{fig4}
\end{figure}

\newpage

Summarizing we can say that there only five digits that remains the same
when we gave them a $180^{o}$ rotation and these are

\begin{figure}[htbp]
\centerline{\includegraphics[width=0.87in,height=0.17in]{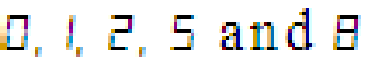}.}
\label{fig5}
\end{figure}

In this paper we shall present different kind of universal magic squares of
order $8\times 8$, $9\times 9$, $16\times 16$ and $25\times 25$ using only the five digits given above.

\bigskip
\subsection{Universal and Bimagic Squares}

\bigskip
Here below are some definitions.

\bigskip
\noindent \textbf{$\bullet$ Magic square}

\bigskip
\noindent  A magic square is a collection of numbers put as a square matrix,
where the sum of elements of each row, sum of elements of each column or sum
of elements of each of two principal diagonals are the same. For
simplicity, let us write this sum as \textbf{S1}.

\bigskip
\noindent \textbf{$\bullet$ Bimagic square}

\bigskip
\noindent Bimagic square is a magic square where the sum of square of each
element of rows, columns or two principal diagonals are the same. For
simplicity, let us write this sum as \textbf{S2}.
\bigskip

\noindent \textbf{$\bullet$ Universal magic square}

\bigskip
\noindent Universal magic square is a magic square with the following properties:

\begin{itemize}
\item[(i)] \textbf{Upside down}, i.e., if we rotate it to $180^{o}$ degrees, it remains magic square again;
\item[(ii)] \textbf{Mirror looking}, i.e., if we put it in front of mirror or see from
the other side of the glass, or see on the other side of the paper, it
always remains the magic square.
\end{itemize}

Here below are the universal bimagic squares of different orders.

\newpage
\noindent \textbf{$\bullet$ Universal bimagic square of order $8 \times 8$ using only the digits 1 and 8}

\begin{figure}[htbp]
\centerline{\includegraphics[width=5.28in,height=2.00in]{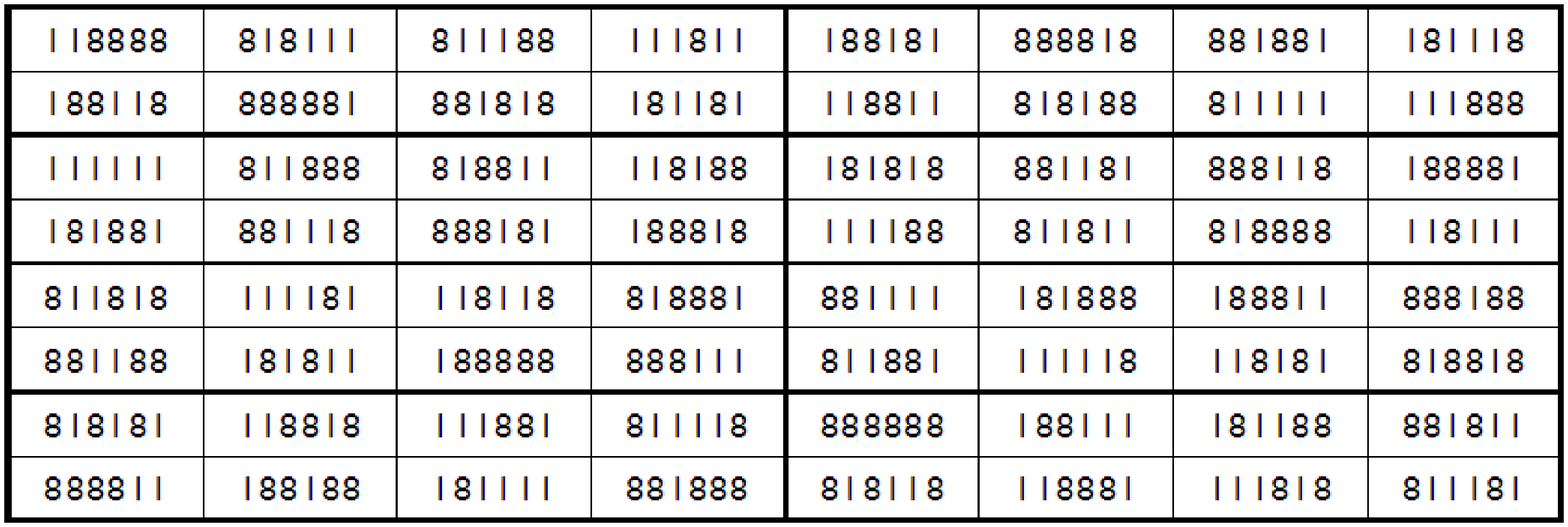}}
\label{fig6}
\end{figure}

\noindent S1:= 3999996\\
S2:=2989894989900

\bigskip
Also we have sum of each block of $2\times 4=3999996$ and the square of sum of each
term in each block of $2 \times 4=2989894989900$

\bigskip
\noindent \textbf{$\bullet$ Universal bimagic square of order $8 \times 8$ using only the digits 2 and 5}

\begin{figure}[htbp]
\centerline{\includegraphics[width=5.28in,height=2.00in]{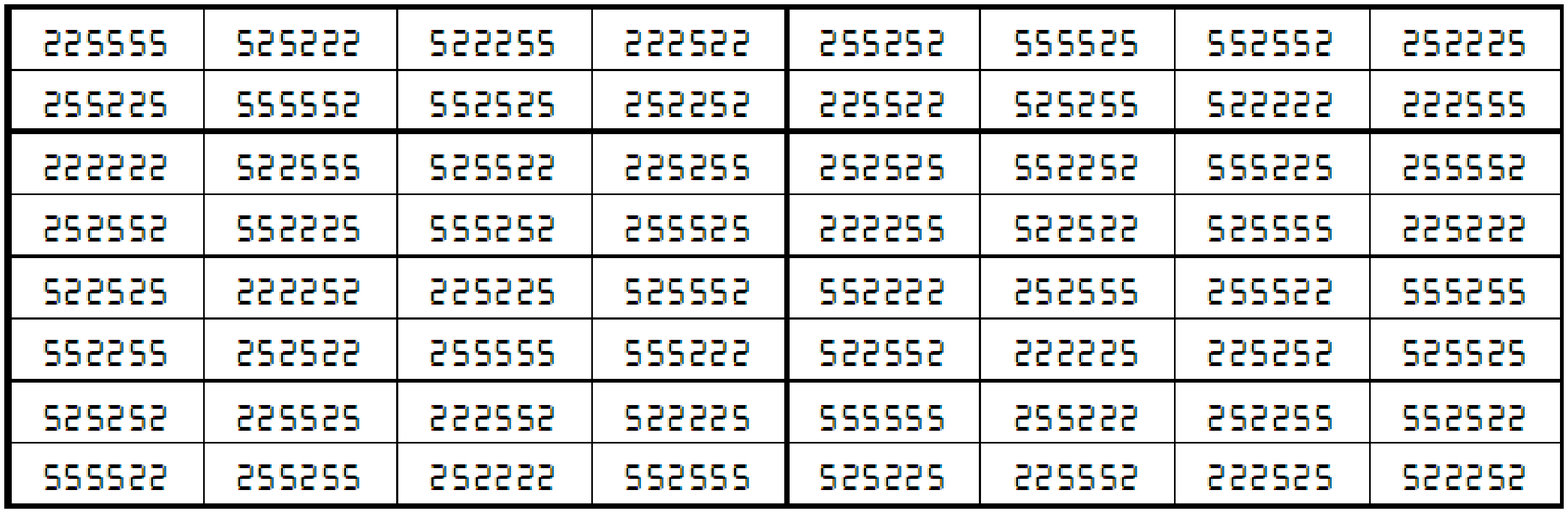}}
\label{fig7}
\end{figure}

\noindent S1:=3111108\\
S2:= 1391692305276

\bigskip
Also we have sum of each block of $2\times 4=3111108$ and the square of sum of each
term in each block of $2 \times 4=1391692305276$

\newpage
\noindent \textbf{$\bullet$ Universal bimagic squares of order $9 \times 9$ using only the digits 1, 2 and 5}

\begin{figure}[htbp]
\centerline{\includegraphics[width=5.28in,height=1.98in]{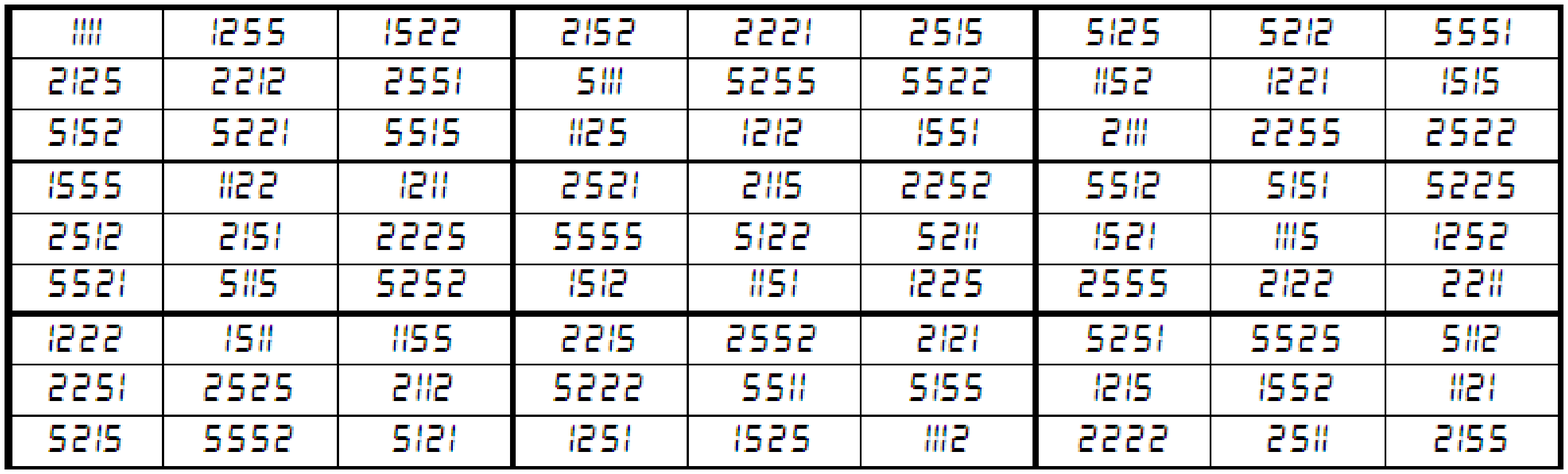}}
\label{fig8}
\end{figure}

\noindent S1:=26664\\
S2:= 105259179

\bigskip
Also we have sum of each block of $3\times 3=26664$ and the square of sum of each
term in each block of $3 \times 3=105259179$

\bigskip
\noindent \textbf{$\bullet$ Universal bimagic square of order $9 \times 9$ using only the digits 2, 5 and 8}

\begin{figure}[htbp]
\centerline{\includegraphics[width=5.28in,height=1.98in]{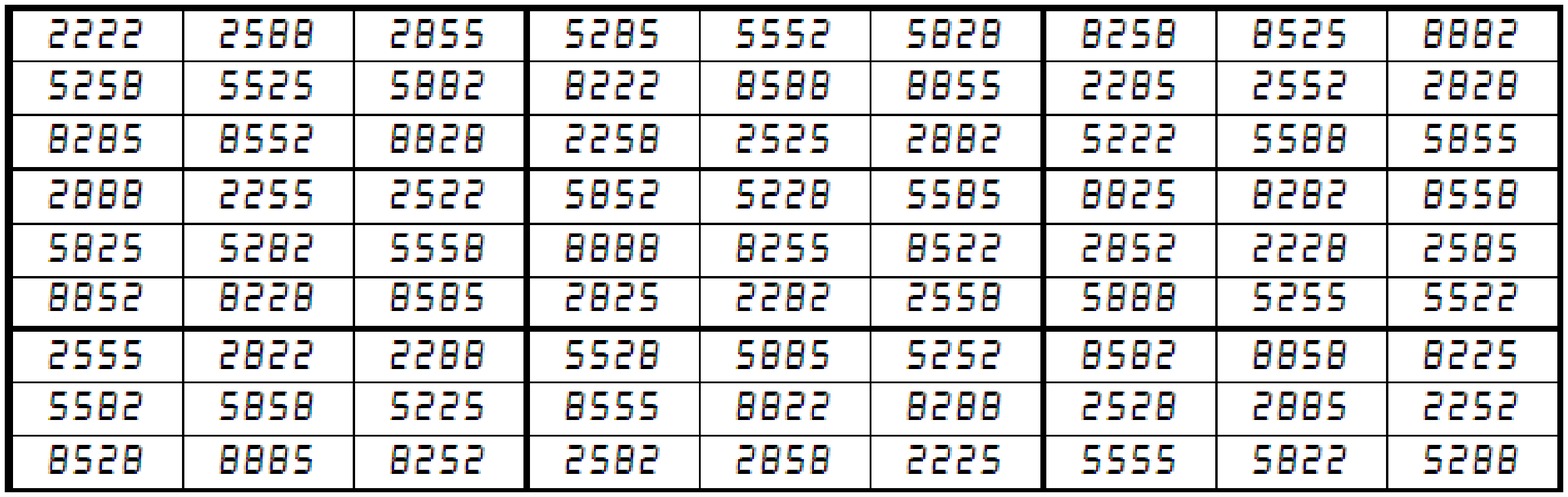}}
\label{fig9}
\end{figure}

\noindent S1:= 49995\\
S2:= 332267679

\bigskip
Also we have sum of each block of $3\times 3=49995$ and the square of sum of each
term in each block of $3 \times 3=332267679$

\newpage
\noindent \textbf{$\bullet$ Universal bimagic square of order $16 \times 16$ using only the digits 1, 2, 5 and 8}

\begin{figure}[htbp]
\centerline{\includegraphics[width=6.28in,height=3.32in]{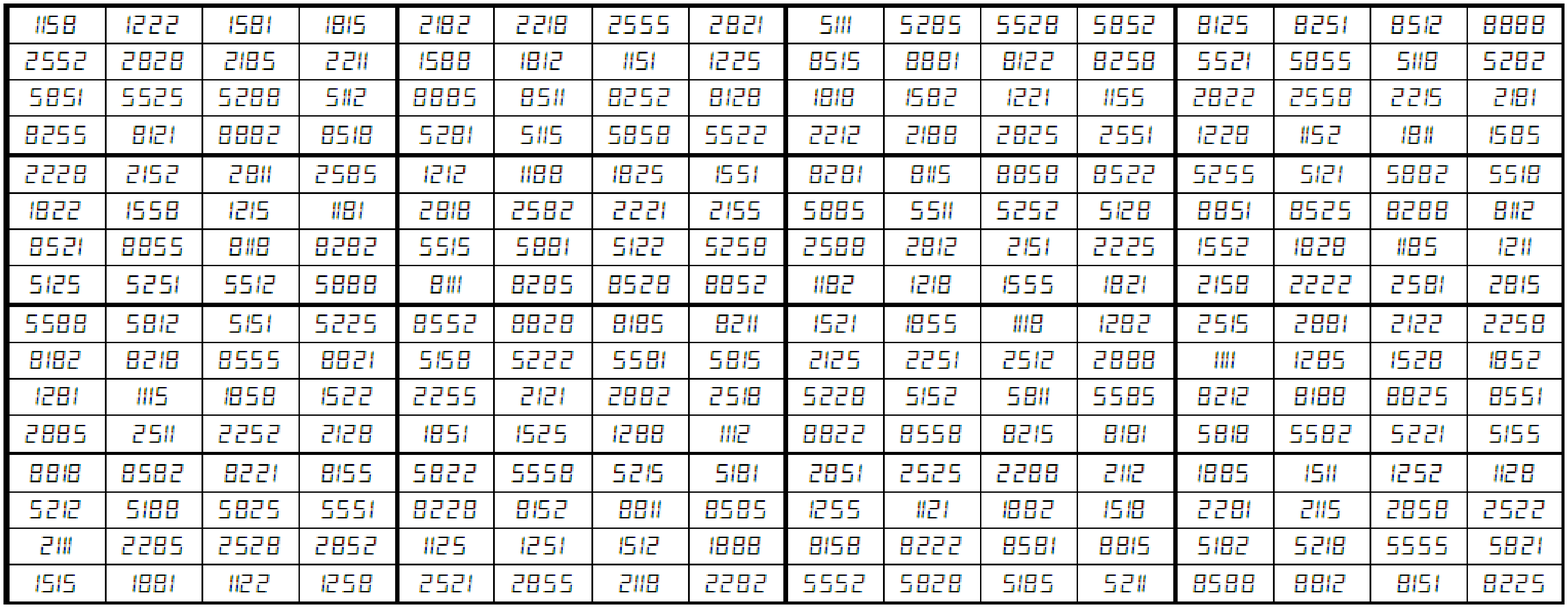}}
\label{fig10}
\end{figure}

\noindent S1:= 71104\\
S2:= 437198296

\bigskip
Also we have sum of each block of $4\times 4= 711045$ and the square of sum of each
term in each block of $4 \times 4=437198296$

\newpage
\noindent \textbf{$\bullet$ Universal bimagic square of order $25 \times 25$ using only the digits 0, 1, 2, 5 and 8}

\begin{figure}[htbp]
\centerline{\includegraphics[width=6.28in,height=3.80in]{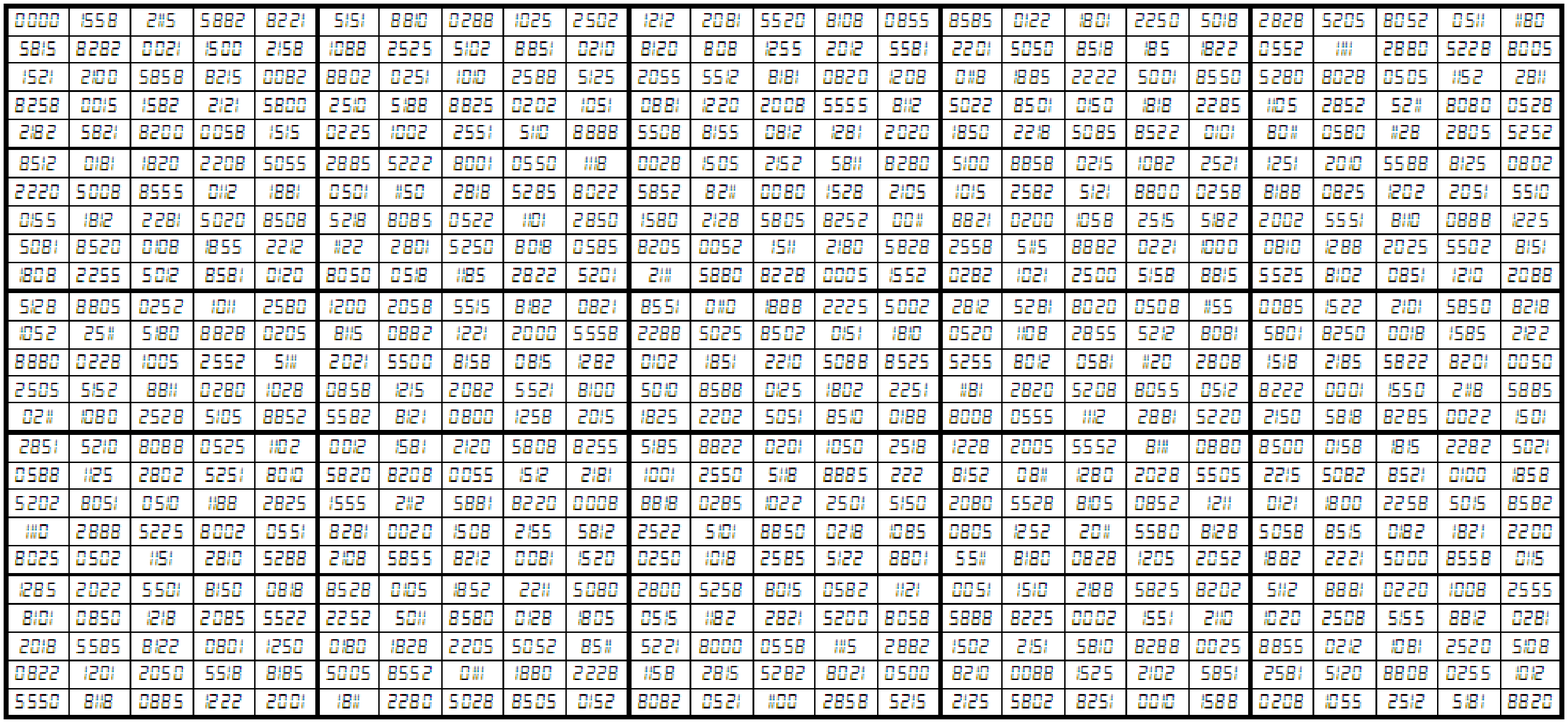}}
\label{fig11}
\end{figure}

\noindent S1:=88880\\
S2:=532147790

\bigskip
Also we have sum of each block of $5\times 5= 88880$ and the square of sum of each
term in each block of $5 \times 5=532147790$

\bigskip
We observe that the above universal bimagic square of order $25 \times 25$ is also is \textit{pendiagonal}. If we write the sum S1 as 088880 then it also becomes upside down and mirror looking. See below:

\newpage

\begin{figure}[htbp]
\centerline{\includegraphics[width=6.28in,height=4.13in]{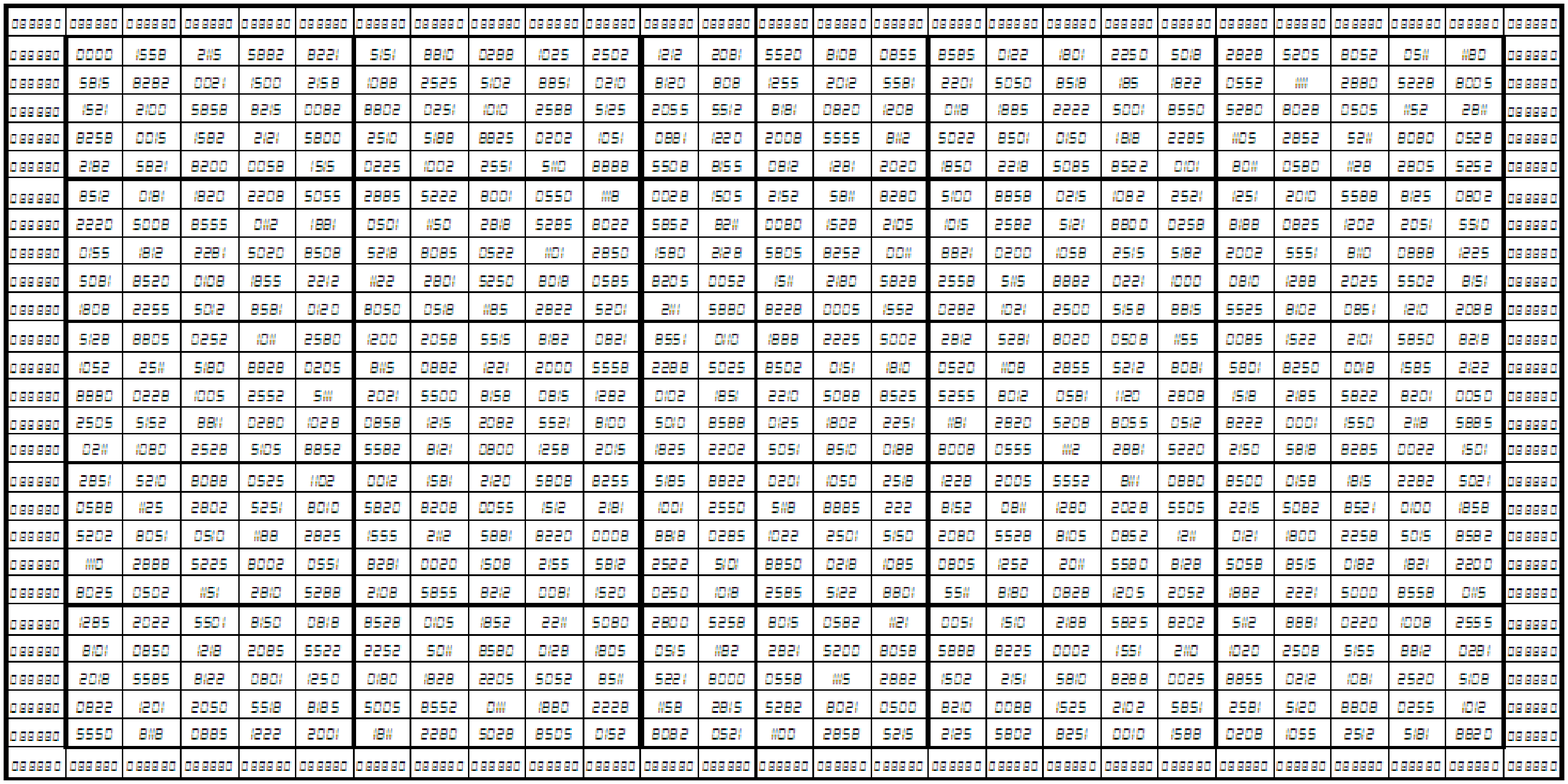}}
\label{fig12}
\end{figure}

We have given above the \textit{universal bimagic squares} using only the five digits 0, 1, 2, 5 and 8. For more examples of \textit{universal bimagic squares} refer to Taneja \cite{tan3}.

\bigskip
If we use the seven digits 0, 1, 2, 5, 6, 8 and 9, then this gives us the conditions to produce a \textit{upside down bimagic square} of order $49\times 49$. In view of 6 and 9, we can't have it as mirror looking, because 6 as well as 9 don't produce any number when we see them in mirror. This study we have presented separately in another work (ref. Taneja \cite{tan4})

\bigskip
For more studies on magic and bimagic squares, we suggest to the readers
the two sites \cite{boyer}, \cite{her} where one can find a good collection of work,
papers, books, etc. The idea of \textit{universal bimagic square} is presented for the first time here.

\begin{center}
---------------------------
\end{center}

\end{document}